\documentclass[11pt]{article}
\usepackage{amsfonts}
\usepackage{latexsym}
\usepackage{cite}
\usepackage{amsmath,amsfonts,latexsym,amssymb,color}
\usepackage[mathscr]{eucal}
\usepackage{cases}
\usepackage{amsthm}

\usepackage[bf,small]{caption2}
\usepackage{float}
\usepackage{graphicx}
\usepackage{amsmath}
\usepackage{amssymb}
\usepackage[all]{xy}

\newtheorem{theorem}{theorem}[section]
\newtheorem{thm}[theorem]{Theorem}

\newtheorem{prob}[theorem]{Problem}
\newtheorem{prop}[theorem]{Proposition}
\newtheorem{cor}[theorem]{Corollary}

\newtheorem{exmp}[theorem]{Example}

\newtheorem{rmk}[theorem]{Remark}
\newtheorem{nota}[theorem]{Notation}

\begin{document}

\title{\textbf{Quotients of polynomial rings and regular t-balanced Cayley maps on abelian groups}}
\author{\Large Haimiao Chen \\
\normalsize \em{Mathematics, Beijing Technology and Business University, Beijing, China} \\
{\em \small chenhm@math.pku.edu.cn}}
\date{}
\maketitle

\begin{abstract}
  Given a finite group $\Gamma$, a regular $t$-balanced Cayley map (RBCM$_{t}$ for short) is a regular Cayley map $\mathcal{CM}(G,\Omega,\rho)$ such that $\rho(\omega)^{-1}=\rho^{t}(\omega)$ for all $\omega\in\Omega$.
  In this paper, we clarify a connection between quotients of polynomial rings and RBCM$_{t}$'s on abelian groups, so as to propose a new approach for classifying RBCM$_{t}$'s. We obtain many new results, in particular, a complete classification for RBCM$_{t}$'s on abelian 2-groups.

\end{abstract}

\section{Introduction}

Suppose $\Gamma$ is a finite group and $\Omega$ is a generating set such that $\Omega$ does not contain the identity $e$, and $\omega^{-1}\in\Omega$ whenever $\omega\in\Omega$. The \emph{Cayley graph} $\textrm{Cay}(\Gamma,\Omega)$ is the graph having the vertex set $\Gamma$ and the arc set $\Gamma\times\Omega$, where $(\eta,\omega)$ denotes the unique arc from $\eta$ to $\eta\omega$.

Let $\rho$ be a cyclic permutation on $\Omega$. It induces a permutation on $\Gamma\times\Omega$ via $(\eta,\omega)\mapsto(\eta,\rho(\omega))$,
and equips each vertex $\eta$ with a cyclic order, which means a cyclic permutation on the set of arcs emanating from $\eta$.
This determines an embedding of ${\rm Cay}(\Gamma,\Omega)$ into a unique closed oriented surface, characterized by the property that each connected component of the complement of the graph is a disk. The embedding, denoted by $\mathcal{CM}(\Gamma,\Omega,\rho)$, is called a {\it Cayley map} (on $\Gamma$).

An \emph{isomorphism} $\mathcal{CM}(\Gamma,\Omega,\rho)\to\mathcal{CM}(\Gamma',\Omega',\rho')$ is an isomorphism between their underlying graphs ${\rm Cay}(\Gamma,\Omega)\to{\rm Cay}(\Gamma',\Omega')$ which can be extended to an orientation-preserving homeomorphism between the supporting surfaces; it can be written as
$\alpha:\Gamma\times\Omega\to\Gamma'\times\Omega'$
such that
\begin{align}
\alpha\circ\rho=\rho'\circ\alpha.   \label{eq:iso}
\end{align}

A Cayley map is said to be \emph{regular} if
for any two arcs, there exists a unique automorphism sending one arc to the other.
It was shown in \cite{JS02} that a Cayley map $\mathcal{CM}(\Gamma,\Omega,\rho)$ is regular if and only if there exists a  bijective function (called {\it skew-morphism}) $\varphi:\Gamma\to\Gamma$ with an associated {\it power function} $\pi:\Gamma\to\{1,\ldots,|\Omega|\}$ such that $\varphi|_{\Omega}=\rho$, $\varphi(e)=e$ and
$$\varphi(\eta\mu)=\varphi(\eta)\varphi^{\pi(\eta)}(\mu) \qquad \text{for\ all} \quad \eta,\mu\in\Gamma.$$
Note that the function
$$\Gamma\times\Omega\to\Gamma\times\Omega,  \qquad (\eta,\omega)\mapsto(\varphi(\eta),\rho^{\pi(\eta)}\omega)$$
defines an automorphism of $\mathcal{CM}(\Gamma,\Omega,\rho)$.

Regular Cayley maps are highly symmetric objects of combinatorial and topological nature, constructed from groups.
It is natural to consider
\begin{prob}
\rm Given a group $\Gamma$, classify all regular Cayley maps on $\Gamma$.
\end{prob}
In general this is a difficult problem. One reason is, as we shall see, the category of regular Cayley maps is very diverse; already RBCM$_{t}$'s on abelian groups form quite a large variety, let alone general regular Cayley maps.
Only when $\Gamma$ is cyclic, a complete classification was obtained recently in \cite{CT14} (2014).

As a first step, one may focus on some special classes.
Let $t$ be an integer. A regular Cayley map $\mathcal{CM}(\Gamma,\Omega,\rho)$ is said to be \emph{$t$-balanced} if
\begin{align}
\rho(\omega)^{-1}=\rho^{t}(\omega^{-1}) \ \ \ \text{for\ all\ \ } \omega\in\Omega;  \label{eq:t-balanced}
\end{align}
in particular, it is called {\it balanced} if $t=1$ and {\it anti-balanced} if $t=-1$.
Note that, not $t$, but its residue modulo $|\Omega|$ is relevant.
From now on we assume $t>0$, and abbreviate ``regular $t$-balanced Cayley map" to ``RBCM$_{t}$". 

Recall some well-known facts on RBCM$_{t}$'s in the following Proposition.
Actually (a) is Theorem 2 of \cite{SS92}; for (b),(c), one may refer to  Lemma 3.1 of \cite{FJW11}, Lemma 2.3 of \cite{KKF06} and Lemma 4.1 of \cite{CJT07-t}.

\begin{prop}  \label{prop:RBCMt}
{\rm(a)} A Cayley map $\mathcal{CM}(\Gamma,\Omega,\rho)$ is a RBCM$_{1}$ if and only if $\rho$ can be extended to an automorphism of $\Gamma$.

{\rm(b)} Suppose $1<t<|\Omega|$. A Cayley map $\mathcal{CM}(\Gamma,\Omega,\rho)$ is a RBCM$_{t}$ if and only if $\pi(\omega)=t$ for all $\omega\in\Omega$ and $\pi(\eta)\in\{1,t\}$ for all $\eta\in\Gamma$.

{\rm(c)} When the conditions in {\rm(b)} are satisfied,
$\Gamma^{+}:=\{\eta\in \Gamma\colon \pi(\eta)=1\}$ is a subgroup of index 2 in $\Gamma$,
consisting of elements which are products of an even number of generators,
$\varphi(\Gamma^{+})=\Gamma^{+}$, and $\varphi|_{\Gamma^{+}}$ is an automorphism.
\end{prop}

Suppose $\mathcal{CM}(\Gamma,\Omega,\rho)$ is a RBCM$_{t}$ with $|\Omega|=m$. Choose an element $\omega\in\Omega$ and denote
$\omega_{i}=\rho^{i}(\omega)$; suppose $\omega^{-1}=\omega_{\ell}$.
By the condition (\ref{eq:t-balanced}),
$$\omega_{i}^{-1}=\omega_{\ell+ti},$$
where the subscript is understood as the residue modulo $m$.
Replacing $i$ by $\ell+ti$, we obtain $\omega_{\ell+ti}^{-1}=\omega_{(t+1)\ell+t^{2}i}$,
hence $\omega_{i}=\omega_{(t+1)\ell+t^{2}i}$.
Since this is true for all $i$, we have
\begin{align}
t^{2}-1\equiv (t+1)\ell\equiv 0\pmod{m}.         \label{eq:t-m-ell}
\end{align}

The integer $\ell$ depends on the choice of $\omega$. If we choose another element $\omega'$, say $\omega'=\omega_{j}$, define $\omega'_{i}=\rho^{i}(\omega')$ and $(\omega')^{-1}=\omega'_{\ell'}$, then $\ell'-\ell=(t-1)j$.
By (\ref{eq:t-m-ell}), $(t-1,m)\mid 2\ell$, (when $t=1$, by $(t-1,m)$ we mean $m$ for convenience), so there are two possibilities.
\begin{enumerate}
  \item [\rm(I)] If $(t-1,m)\nmid\ell$, then $2\mid (t-1,m)$ and we may take $j$ so that $\ell'=(t-1,m)/2$.
  \item [\rm(II)] If $(t-1,m)\mid\ell$, then we may take $j$ so that $\ell'=(t-1,m)$.
\end{enumerate}
We say that a RBCM$_{t}$ satisfying case I or II is of {\it type I} or {\it type II}, respectively.

\begin{rmk}
\rm To classify RBCM$_t$'s up to isomorphism, we may assume $\ell\in\{(t-1,m)/2,(t-1,m)\}$, as done throughout this paper.
\end{rmk}

\begin{rmk}   \label{rmk:iso}
\rm Observe that a RBCM$_{t}$ $\mathcal{CM}(\Gamma,\{\omega_1,\ldots,\omega_m\},\rho)$ has type II if and only if $\omega_i$ is of order 2 for some $i$; especially, a RBCM$_{1}$ on an abelian group $\Gamma$ has type II if and only if $\Gamma$ is an elementary 2-group.
Thus RBCM$_{t}$'s of different types cannot be isomorphic.

According to Lemma 2.4 of \cite{KKF06}, two RBCM$_t$'s $\mathcal{CM}(\Gamma,\{\omega_1,\ldots,\omega_{m}\},\rho)$ and $\mathcal{CM}(\Gamma',\{\omega'_1,\ldots,\omega'_{m}\},\rho')$ of the same type are isomorphic if and only if there exists an isomorphism $\Gamma\to\Gamma'$ sending $\omega_{i}$ to $\omega'_{i}$.

We shall mention a delicate point: it is RBCM$_t$'s for each given $t$ that are classified in this paper;  the possibility that two RBCM$_t$'s for different $t$'s are isomorphic as maps is not excluded.
\end{rmk}

There have been several results on classifying RBCM$_{t}$'s. Wang and Feng \cite{WF05} (2005)  classified RBCM$_{1}$'s for cyclic groups, dihedral groups, and generalized quaternion groups. Conder, Jajcay and Tucker  \cite{CJT07-t} (2007) classified RBCM$_{-1}$'s on abelian groups. Progress on classifying general RBCM$_{t}$'s includes: Kwak, Kwon and Feng \cite{KKF06} (2006) for dihedral groups, Kwak and Oh \cite{KO08} (2008) for dicyclic groups, Oh \cite{Oh09} (2009)  for semi-dihedral groups, and Kwon \cite{Kwo13}  (2013)  for cyclic groups.

In this paper, we aim to classify RBCM$_{t}$'s for abelian groups.
In Section 2, we establish a correspondence between RBCM$_{1}$'s and quotients $\mathbb{Z}_{N}[x]/Q$ ($\mathbb{Z}_{N}=\mathbb{Z}/N\mathbb{Z}$, regarded as a quotient ring of $\mathbb{Z}$), with $Q$ an ideal containing $x^{n}+1$; this is a warm-up for the more general situation, and motivates the investigation of ideals of polynomial rings.
Section 3 includes some preliminaries on the polynomial rings $\mathbb{Z}_{p^{k}}[x]$, with $p$ a prime.
As a result of its own interest, the quotients $\mathbb{Z}_{2^{k}}[x]/Q$ with $x^{n}+1\in Q$ are all determined.
In Section 4,  we establish a 1-1 correspondence between RBCM$_{t}$'s (for $t\ne 1$) on abelian groups and quotients $\mathbb{Z}_{N}[x]/Q$ with $Q$ satisfying specific conditions. This enables us to study RBCM$_{t}$'s using tools from commutative algebra. With modest efforts, we succeed in completely classifying type II RBCM$_{t}$'s on abelian groups, and reduce the classification of type I RBCM$_{t}$'s to the problem of determining quotients $\mathbb{Z}_{p^{k}}[x]/Q$ with $p\ne 2$, which is left as an open problem.
As a special class, RBCM$_{t}$'s on abelian 2-groups are completely classified.

This approach is expected to shed light on structure of more general regular Cayley maps on abelian groups.
By Lemma 5.1 of \cite{CJT07}, if $\mathcal{CM}(\Gamma,\Omega,\rho)$ is a regular Cayley map with $\Gamma$ abelian and $\varphi$ is the associated skew-morphism, then $\ker\pi=\{\eta\in\Gamma\colon \pi(\eta)=1\}$ is a subgroup such that $\varphi(\ker\pi)=\ker\pi$ and $\varphi|_{\ker\pi}$ is an automorphism. Further research in this direction may reveal the structure of $\ker\pi$, possibly recognizing it as a quotient of some commutative ring of certain kind.

\medskip

\textbf{Notational convention} \

For $N>1$, via the quotient $\mathbb{Z}\twoheadrightarrow\mathbb{Z}_{N}$, an integer $u$ is also considered as an element of $\mathbb{Z}_{N}$.

For an ideal $Q\leq\mathbb{Z}_N[x]$, denote $\mathbb{Z}_N[x]/Q$ by $F_N(Q)$.

Given $\alpha\in\mathbb{Z}_{N}[x]$. The image of $\alpha$ under the quotient $\mathbb{Z}_{N}[x]\twoheadrightarrow F_{N}(Q)$ is denoted by $\alpha+Q$, and shortend as $\overline{\alpha}$ whenever there is no ambiguity;
the bar is omitted if $\alpha\in\mathbb{Z}_{N}$.
For $N'\mid N$, the image of $\alpha$ under the reduction $\mathbb{Z}_{N}[x]\twoheadrightarrow\mathbb{Z}_{N'}[x]$  is denoted by $\mathcal{R}_{N'}(\alpha)$.

For an integer $n>1$, let $\mathcal{I}(N,n)$ denote the set of ideals $Q\le\mathbb{Z}_{N}[x]$ such that
$x^{n}+1\in Q$ and $x^{u}+1\notin Q$ whenever $0<u<n$.

For a prime $p$ and $u\ne 0$, let $d_{p}(u)$ denote the largest $s$ with $p^{s}\mid n$.

If $\alpha$ is an element of some ring and $u>1$, denote $[u]_{\alpha}=1+\alpha+\cdots+\alpha^{u-1}$.

For an abelian group $\Gamma$, denote its identity by 0, and denote its exponent
(the least common multiple of the orders of elements) by $\exp(\Gamma)$.

A Cayley map $\mathcal{CM}(\Gamma,\Omega,\rho)$ is shortend as $\mathcal{CM}(\Gamma,\Omega)$ if $\Omega$ can be written as
$\{\omega_1,\ldots,\omega_m\}$ and $\rho(\omega_i)=\omega_{i+1}$ for all $i$.

\section{Regular balanced Cayley maps on abelian groups and quotients of polynomial rings}

Given $Q\in\mathcal{I}(N,n)$, define the subset $\Omega_{Q}\subset F_{N}(Q)$ by
\begin{align}
\Omega_{Q}=\left\{\begin{array}{ll}
\{1,\overline{x},\ldots,\overline{x}^{2n-1}\},  & N>2, \\
\{1,\overline{x},\ldots,\overline{x}^{n-1}\},   & N=2.
\end{array}\right.
\end{align}
Let $\mathcal{M}_{Q}=\mathcal{CM}(F_{N}(Q),\Omega_{Q})$. When $N>2$, by the definition of $\mathcal{I}(N,n)$, the $\overline{x}^i$'s in $\Omega_Q$ are distinct from each other, and $\overline{x}^{i+n}=-\overline{x}^i$ for each $i$, so $\mathcal{M}_Q$ is a $2n$-valent type I RBCM$_{1}$. Similarly, when $N=2$, $\mathcal{M}_Q$ is an $n$-valent type II RBCM$_{1}$.

If $\mathcal{M}_{Q}\cong\mathcal{M}_{Q'}$ with $Q,Q'\in\mathcal{I}(N,n)$,
then by Remark \ref{rmk:iso} there exists an isomorphism
$F_{N}(Q)\to F_{N}(Q')$ sending $x^{i-1}+Q$ to $x^{i-1}+Q'$ for each $i$, hence $Q=Q'$.

\begin{thm} \label{thm:standard-RBCM}
When $N>2$ (resp. $N=2$), each $2n$-valent type I (resp. $n$-valent type II) RBCM$_{1}$ on an abelian group of exponent $N$ is isomorphic to $\mathcal{M}_{Q}$ for a unique  $Q\in\mathcal{I}(N,n)$.
\end{thm}

\begin{proof}
Suppose $N>2$ and $\mathcal{CM}(\Gamma,\Omega,\rho)$ is a $2n$-valent type I RBCM$_{1}$ with
$\Omega=\{\pm\omega_{1},\ldots,\pm\omega_n\}$, $\rho(\pm\omega_{i})=\pm\omega_{i+1}$ for $1\le i<n$ and $\rho(\pm\omega_{n})=\mp\omega_{n}$.

As an abelian group, $F_{N}((x^{n}+1))$ is generated by $\overline{x}^{i}$ ($0\le i<n$),
and each element can be uniquely written as $\sum\limits_{i=1}^{n}a_{i}\overline{x}^{i-1}$ with $a_{i}\in\mathbb{Z}_{N}$.
Set
$$\overline{Q}=\left\{\sum\limits_{i=1}^{n}a_{i}\overline{x}^{i-1}\colon a_{1},\ldots,a_{n}\in\mathbb{Z}_{N}, \ \sum\limits_{i=1}^{n}a_{i}\omega_{i}=0\right\};$$
the scalar products $a_{i}\omega_{i}$ are well-defined, as $\exp(\Gamma)=N$.
Then $\overline{Q}$ is an ideal of $F_{N}((x^{n}+1))$: by Proposition \ref{prop:RBCMt} (a),  $\rho$ can be extended to an automorphism of $\Gamma$ which we denote by $\phi$;
for each $\alpha=\sum\limits_{i=1}^{n}a_{i}\overline{x}^{i-1}\in\overline{Q}$,
since
$$-a_{n}\omega_{1}+a_{1}\omega_{2}+\cdots+a_{n-1}\omega_{n}=\sum\limits_{i=1}^{n}a_{i}\rho(\omega_{i})
=\phi\left(\sum\limits_{i=1}^{n}a_{i}\omega_{i}\right)=\phi(0)=0,$$
we have $\overline{x}\alpha=-a_{n}+\sum\limits_{i=1}^{n-1}a_{i}\overline{x}^{i}\in\overline{Q}$.

Let $Q=\{\alpha\in\mathbb{Z}_{N}[x]\colon\overline{\alpha}\in\overline{Q}\}$; it is an ideal of $\mathbb{Z}_{N}[x]$ containing $x^{n}+1$, and does not contain $x^{u}+1$ for any $u<n$, since $\omega_{u+1}+\omega_{1}\neq 0$. Hence $Q\in\mathcal{I}(N,n)$.
We now have a group isomorphism
$$\Gamma\to F_{N}(Q), \qquad \omega_{i}\mapsto \overline{x}^{i-1},$$
which defines an isomorphism $\mathcal{CM}(\Gamma,\Omega,\rho)\cong\mathcal{M}_{Q}.$

The case $N=2$ is similar.
\end{proof}

\section{Ideals of the ring $\mathbb{Z}_{p^{k}}[x]$}

For each positive integer $d$, let
$$\Psi_{d}(x)=\prod\limits_{1\leq j\le d\colon (j,d)=1}\left(x-(e^{\frac{2\pi i}{d}})^{j}\right).$$
It is well-known that (see \cite{Mor96} Lemma 7.6) $\Psi_{d}(x)\in\mathbb{Z}[x]$, and for any $n$,
$$x^{n}-1=\prod\limits_{d|n}\Psi_{d}(x).$$

\begin{nota}
\rm For a positive integer $d$ and a prime $p$ with $p\nmid d$, let $o_{p}(d)$ denote the order of $p$ in the multiplicative group $\mathbb{Z}_{d}^{\times}$ of units in $\mathbb{Z}_{d}$; choose and fix a set $\Xi_{p}(d)$ of coset representatives of the subgroup
$\langle p\rangle\leqslant\mathbb{Z}_{d}^{\times}$.
\end{nota}

For the following fact, one can refer to 
\cite{LN83} Theorem 2.45 and Theorem 2.47. For $p\nmid d$, over the algebraic closure of $\mathbb{Z}_{p}$, $\Psi_{d}(x)$ is factorized as
$$\Psi_{d}(x)=\prod\limits_{1\leq i<d\colon (i,d)=1}(x-\xi_{d}^{i}),$$
where $\xi_{d}$ is a $d$-th primitive root of unity, and over $\mathbb{Z}_{p}$,
$\Psi_{d}(x)$ is decomposed into $|\mathbb{Z}_{d}^{\times}|/o_{p}(d)$ monic irreducible polynomials of degree $o_{p}(d)$:
$$\Psi_{d}(x)=\prod\limits_{\ell\in\Xi_{p}(d)}\psi_{d,\ell}, \hspace{5mm}   \text{with\ \ }
\psi_{d,\ell}=\prod\limits_{j=1}^{o_{p}(d)}(x-\xi_{d}^{\ell p^{j}})\in\mathbb{Z}_{p}[x].$$

Thus, for a positive integer $n'$ coprime to $p$, in $\mathbb{Z}_{p}[x]$ we have
\begin{align}
x^{n'}-1=\prod\limits_{d\mid n'}\prod\limits_{\ell\in\Xi_{p}(d)}\psi_{d,\ell}, \label{eq:fac1}
\end{align}
and for $n=p^{r}n'$ with $p\nmid n'$,
\begin{align}
x^{n}-1=(x^{n'}-1)^{p^{r}}=\prod\limits_{d\mid n'}\prod\limits_{\ell\in\Xi_{p}(d)}\psi_{d,\ell}^{p^{r}}.   \label{eq:fac2}
\end{align}

\begin{rmk}  \label{rmk:divide}
\rm It is known that $\psi_{d_{1},\ell_{1}}\neq\psi_{d_{2},\ell_{2}}$ unless $d_{1}=d_{2}$ and $\ell_{1}=\ell_{2}$, so for any $u>0$, $\psi_{d,\ell}^{u}\mid x^{n}-1$ if and only if $u\le p^{r}$.
\end{rmk}

\begin{nota}
\rm  For $n=p^{r}n'$ with $p\nmid n'$, let $\Lambda(p,n)=\bigcup\limits_{d\mid n'}\Xi_{p}(d).$
Use $\lambda$ to denote a general element of $\Lambda(p,n)$;
if $\lambda=(d,\ell)$, then we denote $\psi_{\lambda}=\psi_{d,\ell}$.
Let $\o$ denote the special value $(1,1)$ so that $\psi_{\o}=x-1$.
\end{nota}

Now consider the factorizations in $\mathbb{Z}_{p^{k}}[x]$.
According to \cite{GG99} Algorithm 15.17, Theorem 15.18 (see also Algorithm 15.10, Theorem 15.12 and Theorem 15.14), given $\alpha\in\mathbb{Z}[x]$ and $\alpha_{1},\ldots,\alpha_{s}\in\mathbb{Z}_{p}[x]$ such that $(\alpha_{i},\alpha_{j})=1$ for all $i\ne j$ and $\alpha=\alpha_{1}\cdots\alpha_{s}$ in $\mathbb{Z}_{p}[x]$, there exist unique $\tilde{\alpha}_{1},\ldots,\tilde{\alpha}_{s}\in\mathbb{Z}_{p^{k}}[x]$ such that $\mathcal{R}_{p}(\tilde{\alpha}_{j})=\alpha_{j}, j=1,\ldots,s$, and $\alpha=\tilde{\alpha}_{1}\cdots\tilde{\alpha}_{s}$ in $\mathbb{Z}_{p^{k}}[x]$; the ``lifts" $\tilde{\alpha}_{j}$ can be found by a practical method called ``Hensel lifting". One may also refer to \cite{BF02, McG01}.

As special cases and their simple consequences, we have
\begin{prop}  \label{prop:fundamental}
{\rm(a)} For each $\lambda\in\Lambda(p,n)$, in $\mathbb{Z}_{p^{k}}[x]$ there exists a unique monic factor
$\psi_{\lambda,r}$ of $x^{n}-1$ with $\mathcal{R}_{p}(\psi_{\lambda,r})=\psi_{\lambda}^{p^{r}}$.

{\rm(b)} We have $\psi_{\lambda,r}\mid\psi_{\lambda,r+1}$ and $\phi_{\lambda,r}:=\psi_{\lambda,r+1}/\psi_{\lambda,r}$ is the unique monic factor of $[p]_{x^{n'p^{r}}}$ such that $\mathcal{R}_{p}(\phi_{\lambda,r})=(\psi_{\lambda})^{p^{r}(p-1)}$.

{\rm(c)} Denoting
$\overline{\Lambda}(p,n)=\Lambda(p,2n)-\Lambda(p,n)$, we have
\begin{align}
x^{n}-1=\prod\limits_{\lambda\in\Lambda(p,n)}\psi_{\lambda,r},  \hspace{10mm}
x^{n}+1=\left\{\begin{array}{ll}
\prod\limits_{\lambda\in\overline{\Lambda}(p,n)}\psi_{\lambda,r}, & p\ne 2, \\
\prod\limits_{\lambda\in\Lambda(2,n)}\phi_{\lambda,r}, & p=2.
\end{array}\right.
\end{align}
\end{prop}

\begin{exmp}   \label{exmp}
\rm As a special case, if $p=2$, then
$$\psi_{\o,r}=x^{2^{r}}-1, \hspace{10mm}   \phi_{\o,r}=x^{2^{r}}+1.$$
\end{exmp}

\begin{prop}  \label{prop:iso-of-ring}
{\rm(a)} Let $Q\le\mathbb{Z}_{p^{k}}[x]$ be an ideal containing $x^{n}-1$. For each $\lambda\in\Lambda(p,n)$, let $Q(\lambda)=Q+(\psi_{\lambda,r})$.   There is a ring isomorphism
\begin{align*}
F_{p^{k}}(Q)\cong\prod\limits_{\lambda\in\Lambda(p,n)}F_{p^{k}}(Q(\lambda)), \qquad
\overline{\alpha}\mapsto\{\alpha+Q(\lambda)\}_{\lambda\in\Lambda(p,n)}.
\end{align*}

{\rm (b)} Suppose $p\ne 2$ and $Q\le\mathbb{Z}_{p^{k}}[x]$ is an ideal containing $x^{n}+1$. For each $\lambda\in\overline{\Lambda}(p,n)$, let $Q(\lambda)=Q+(\psi_{\lambda,r})$. There is a ring isomorphism
\begin{align*}
F_{p^{k}}(Q)\cong\prod\limits_{\lambda\in\overline{\Lambda}(p,n)}F_{p^{k}}(Q(\lambda)), \qquad
\overline{\alpha}\mapsto\{\alpha+Q(\lambda)\}_{\lambda\in\overline{\Lambda}(p,n)}.
\end{align*}

{\rm (c)} Suppose $Q\le\mathbb{Z}_{2^{k}}[x]$ is an ideal containing $x^{n}+1$. For each $\lambda\in\Lambda(2,n)$, let $Q(\lambda)=Q+(\phi_{\lambda,r})$. There is a ring isomorphism
\begin{align*}
F_{2^{k}}(Q)\cong\prod\limits_{\lambda\in\Lambda(2,n)}F_{2^{k}}(Q(\lambda)), \qquad
\overline{\alpha}\mapsto\{\alpha+Q(\lambda)\}_{\lambda\in\Lambda(2,n)}.
\end{align*}
\end{prop}

\begin{proof}
We only prove (a); the proofs for (b), (c) are similar.

Let $\overline{Q(\lambda)}$ denote the image of $Q(\lambda)$ under the quotient
$\mathbb{Z}_{p^{k}}[x]\twoheadrightarrow F_{p^{k}}(Q)$, so that $\overline{Q(\lambda)}$ is an ideal of $F_{p^{k}}(Q)$.
Whenever $\lambda\neq\lambda'$, $\psi_{\lambda}$ is coprime to $\psi_{\lambda'}$, hence $\psi_{\lambda,r}$ is coprime to $\psi_{\lambda',r}$, implying that $\overline{Q(\lambda)}$ is coprime to $\overline{Q(\lambda')}$. Also, by Proposition \ref{prop:fundamental} (c), $\prod\limits_{\lambda\in\Lambda(p,n)}\overline{Q(\lambda)}=0$ in $F_{p^{k}}(Q)$.
The result then follows from Proposition 1.10 of \cite{AM}.
\end{proof}

Therefore, when $p\ne 2$ (resp. $p=2$), the study of quotients $F_{p^{k}}(Q)$ with $x^{n}+1\in Q$ is reduced to that of quotients $F_{p^{k}}(Q)$ with $\psi_{\lambda,r}\in Q$ (resp. $\phi_{\lambda,r}\in Q$).

\begin{thm}   \label{thm:ideal}
For each $\lambda\in\Lambda(2,n)$, there are exactly $k(2^{r}+1)$ nontrivial ideals of $\mathbb{Z}_{2^{k}}[x]$ that contain $\phi_{\lambda,r}$, namely,
\begin{align}
M_{u,v}^{\lambda}:=(\phi_{\lambda,r},2^{u}\psi_{\lambda,0}^{v}, 2^{u+1}), \quad  u\in\{0,\ldots,k-1\}, \ v\in\{0,\ldots,2^{r}\}.
\end{align}
\end{thm}

\begin{proof}
The result is clear if $k=1$, so we assume $k\ge 2$.
Let $M<\mathbb{Z}_{2^{k}}[x]$ be an ideal containing $\phi_{\lambda,r}$.
Abbreviate $\psi_{\lambda}$, $\psi_{\lambda,0}$, $\phi_{\lambda,r}$ to $\psi$, $\tilde{\psi}$, $\phi$, respectively.

For $w\in\{0,\ldots,k\}$, let
$$M_{w}=\{\alpha\in\mathbb{Z}_{2^{k}}[x]\colon 2^{w}\alpha\in M\}.$$
Since $\mathcal{R}_{2}(M_{w})$ contains $\mathcal{R}_{2}(\phi)=\psi^{2^{r}}$,
we have $\mathcal{R}_{2}(M_{w})=(\psi^{v_{w}})$ for some $v_{w}\leq 2^{r}$;
arbitrarily take $\gamma_{w}\in M_{w}$ with $\mathcal{R}_{2}(\gamma_{w})=\psi^{v_{w}}$
(when $v_{w}=2^{r}$, $\gamma_{w}$ can be taken to be $\phi$),
then $\gamma_{w}=\tilde{\psi}^{v_{w}}+2\alpha_{w}$ for some $\alpha_{w}$.

Note that $v_{0}\geq\cdots\geq v_{k-1}\geq v_{k}=0$, since
\begin{align*}
M=M_{0}\leq M_{1}\leq\cdots\leq M_{k-1}\leq M_{k}=\mathbb{Z}_{2^{k}}[x];
\end{align*}
also, $v_{0}>0$, since $M_{0}=M\neq\mathbb{Z}_{2^{k}}[x]$.

For each $w\le k-1$, $2M_{w+1}\leq M_{w}$, and the kernel of the composite
$$M_{w}\hookrightarrow\mathbb{Z}_{2^{k}}[x]\stackrel{\mathcal{R}_{2}}\longrightarrow\mathbb{Z}_{2}[x]$$
is $2M_{w+1}$, hence $M_{w}/2M_{w+1}\cong\mathcal{R}_{2}(M_{w})$, and $M_{w}=(\gamma_{w})+2M_{w+1}$.
Thus
\begin{align}
M=(\gamma_{0})+2M_{1}=(\gamma_{0})+(2\gamma_{1})+2^{2}M_{2}=\cdots=(\gamma_{0},\ldots,2^{k-1}\gamma_{k-1}).   \label{eq:M}
\end{align}

If $v_{k-1}=2^{r}$, then $M=(\phi)=(\phi,2^{k-1}\tilde{\psi}^{2^{r}},2^{k})$.

Now suppose $v_{k-1}<2^{r}$. Let $u=\min\{w\colon 0<v_{w}<2^{r}\}$ and let $v=v_{u}$; let $2\beta$ denote the remainder of $\phi$ divided by $\gamma_{u}$. Note that
$$x^{n}+1=(x^{n}-1)+2\equiv 2\pmod{(\tilde{\psi},4)},$$
so that $\phi\not\equiv 0\pmod{(\tilde{\psi},4)}$, hence $(\mathcal{R}_{2}(\beta),\psi)=1$.
This shows $2\in (2\beta,\gamma_{u})\le M_{u}$, so actually $M_{u}=(\tilde{\psi}^{v},2)$.
By (\ref{eq:M}) we have $M=(\psi_{r},2^{u}\tilde{\psi}^{v},2^{u+1})$.
\end{proof}

\begin{rmk}
\rm Suppose $\lambda=(d,\ell)\in\Lambda(2,n)$ and $o_{2}(d)=o$.
Since $\deg(\phi_{\lambda,r})=2^{r}o$ and $\deg(\psi_{\lambda,0}^{v})=vo\leq 2^{r}$, we can take
$$\{\overline{x}^{i}\colon 0\le i<vo\}\cup\{\overline{x}^{j}\overline{\psi_{\lambda,0}}^{v}\colon 0\le j<(2^{r}-v)o\}$$
as a generating set of $\mathbb{Z}_{2^{k}}[x]/M^{\lambda}_{u,v}$; the relations are
\begin{align*}
2^{{u+1}}\overline{x}^{i}=0, \ \ (0\le i<vo ), \ \ \ \text{and\ \ \ }
2^{{u}}\overline{x}^{j}\overline{\psi_{\lambda,0}}^{v}=0, \ \ ( 0\le j<(2^{r}-v)o ).
\end{align*}
Thus, as an abelian group,
\begin{align}
F_{2^{k}}(M^{\lambda}_{u,v})\cong(\mathbb{Z}_{2^{u+1}})^{vo}\times(\mathbb{Z}_{2^{u}})^{(2^{r}-v)o}.
\end{align}

Thanks to Proposition \ref{prop:iso-of-ring} (c), in this way we have found all the possible quotients $F_{2^k}(Q)$ for $Q\le\mathbb{Z}_{2^k}[x]$ containing $x^n+1$.
\end{rmk}

\begin{prob}
\rm For $p\ne 2$ and $\lambda\in\Lambda(p,n)$, determine ideals of $\mathbb{Z}_{p^{k}}[x]$ that contain $\psi_{\lambda,r}$.
\end{prob}

As we shall see in the next section, the answer to this problem will lead to a complete classification of RBCM$_{t}$'s on abelian groups.

\section{Classifying regular $t$-balanced  Cayley maps on abelian groups (with $t\ne 1$)}

\subsection{Standard form}

Given positive integers $N,m,t,\ell$ with $\ell\in\{(t-1,m)/2,(t-1,m)\}$ and $m\mid t^{2}-1$, let $\mathcal{I}(N,\ell,m,t)$ denote the set of ideals $Q\le\mathbb{Z}_{N}[x]$ such that the following holds in $F_{N}[x](Q)$:
\begin{align}
[m]_{\overline{x}}&=0,  \label{eq:x0}  \\
[u]_{\overline{x}}=0&\Leftrightarrow m\mid u,    \label{eq:distinct}   \\
\overline{x}^{(t-1,m)}&=1,  \label{eq:x1} \\
\overline{x}^{\ell}+[t]_{\overline{x}}&=0;    \label{eq:x2}
\end{align}
note that the last two conditions imply
\begin{align}
(\overline{x}^{\ell}+1)(\overline{x}-1)=0.   \label{eq:csq}
\end{align}

Given $Q\in\mathcal{I}(N,\ell,m,t)$, put
\begin{align}
\Gamma_{Q}=\langle F_{N}(Q),\omega\mid 2\omega=-[\ell]_{\overline{x}}\rangle,
\end{align}
i.e., the abelian group generated by $F_{N}(Q)$ and $\omega$ with an additional relation $2\omega=-[\ell]_{\overline{x}}$,
and define
\begin{align*}
&\pi:\Gamma_{Q}\to\{1,t\},  \hspace{10mm} \overline{\alpha}\mapsto 1, \ \ \ \ \ \ \overline{\alpha}+\omega\mapsto t; \\
&\varphi:\Gamma_{Q}\to\Gamma_{Q}, \hspace{13mm} \overline{\alpha}\mapsto\overline{x}\overline{\alpha}, \ \ \ \ \overline{\alpha}+\omega\mapsto\overline{x}\overline{\alpha}+1+\omega;
\end{align*}
in particular, $\varphi([i]_{\overline{x}}+\omega)=[i+1]_{\overline{x}}+\omega$.

One can verify that $\varphi(\mu+\eta)=\varphi(\mu)+\varphi^{\pi(\mu)}(\eta)$ for all $\mu,\eta\in\Gamma$.
Indeed, this is obvious for $\mu=\overline{\alpha}$;
if $\mu=\overline{\alpha}+\omega, \eta=\overline{\beta}$, then
$$\varphi(\mu+\eta)=\overline{x}(\overline{\alpha}+\overline{\beta})+1+\omega=
\overline{x}\overline{\alpha}+1+\omega+\overline{x}^{t}\overline{\beta}=\varphi(\mu)+\varphi^{t}(\eta);$$
if $\mu=\overline{\alpha}+\omega, \eta=\overline{\beta}+\omega$, then
$\varphi^{\pi(\mu)}(\eta)=\overline{x}^{t}\overline{\beta}+[t]_{\overline{x}}+\omega=\overline{x}\overline{\beta}-\overline{x}^{\ell}+\omega$, hence
$$\varphi(\mu)+\varphi^{\pi(\mu)}(\eta)
=\overline{x}\overline{\alpha}+1+\omega+\overline{x}\overline{\beta}-\overline{x}^{\ell}+\omega
=\overline{x}(\overline{\alpha}+\overline{\beta}-[\ell]_{\overline{x}})=\varphi(\mu+\eta).$$
Thus
$\mathcal{M}_{Q}:=\mathcal{CM}(\Gamma_{Q},\{[i]_{\overline{x}}+\omega\colon 1\le i\le m\})$
is an $m$-valent RBCM$_{t}$.

If $\mathcal{M}_{Q}\cong\mathcal{M}_{Q'}$, then similarly as in Section 2, $Q=Q'$.

\begin{thm}
If $\mathcal{M}=\mathcal{CM}(\Gamma,\{\omega_1,\ldots,\omega_{m}\})$ is a RBCM$_{t}$ on an abelian group $\Gamma$ such that $\exp(\Gamma^{+})=N$,  $\rho(\omega_{i})=\omega_{i+1}$ and $\omega_{\ell+ti}=\omega_{i}^{-1}$, $1\le i\le m$, then
$\mathcal{M}\cong\mathcal{M}_{Q}$ for a unique ideal $Q\in\mathcal{I}(N,\ell,m,t)$.
\end{thm}

\begin{proof}
Let $\theta_{i}=\omega_{i}-\omega_{i-1}$, so that $\omega_{i}=\omega_{m}+\sum\limits_{j=1}^{i}\theta_{j}$ for each $i$, and
\begin{align}
\sum\limits_{i=1}^{m}\theta_{i}=0.             \label{eq:sum-of-theta}
\end{align}
Let $\varphi$ denote the skew-morphism, then
\begin{align*}
\varphi(\theta_i)&=\varphi(\omega_i-\omega_{i-1})=\varphi(\omega_i)+\varphi^t(-\omega_{i-1}) \\
&=\varphi(\omega_i)-\varphi(\omega_{i-1})
=\rho(\omega_i)-\rho(\omega_{i-1})=\omega_{i+1}-\omega_i=\theta_{i+1},
\end{align*}
where in the third equality we use (\ref{eq:t-balanced}).
Similarly as in the proof of Theorem \ref{thm:standard-RBCM}, we can show that
$$\overline{Q}:=\left\{\sum\limits_{i=1}^{m}a_{i}\overline{x}^{i-1}\colon a_{i}\in\mathbb{Z}_{N}, \ \sum\limits_{i=1}^{m}a_{i}\theta_{i}=0\right\}\subset F_{N}(([m]_{x}))$$
is an ideal: for each $\overline{\alpha}=\sum\limits_{i=1}^{m}a_{i}\overline{x}^{i-1}\in\overline{Q}$, since
$$\sum\limits_{i=1}^m a_i\theta_{i+1}=\sum\limits_{i=1}^m a_i\varphi(\theta_i)=\varphi\left(\sum\limits_{i=1}^m a_i\theta_i\right)=\varphi(0)=0,$$
(we can pull $\varphi$ outside the summation by Proposition \ref{prop:RBCMt} (c)), we have $\overline{x}\overline{\alpha}\in \overline{Q}$, (note that $\overline{x}^m-1=(\overline{x}-1)[m]_{\overline{x}}=0$).

Thus $Q:=\{\alpha\in\mathbb{Z}_{N}[x]\colon\overline{\alpha}\in\overline{Q}\}$ is an ideal of
$\mathbb{Z}_{N}[x]$ containing $[m]_{x}$, and there is an isomorphism (of abelian groups)
$$\Gamma^{+}\cong F_{N}(Q), \qquad \theta_{i}\mapsto \overline{x}^{i-1}.$$

Written as equations in $F_{N}(Q)$, (\ref{eq:sum-of-theta}) becomes (\ref{eq:x0}).
Since the $\omega_{i}$'s are distinct from each other, we have (\ref{eq:distinct}).
The condition $\omega_{\ell+ti}=\omega_{i}^{-1}$ becomes
\begin{align*}
[\ell+ti]_{\overline{x}}+[i]_{\overline{x}}+2\omega_{m}=0,  \hspace{10mm}  i=1,\ldots,m,
\end{align*}
which is equivalent to $[\ell]_{\overline{x}}=-2\omega_{m}$ and
\begin{align}
[\ell+t(i+1)]_{\overline{x}}+[i+1]_{\overline{x}}=[\ell+ti]_{\overline{x}}+[i]_{\overline{x}}, \ \ \ \ \ i=1,\ldots,m;  \label{eq:cons}
\end{align}
the later can be rewritten as
$$\overline{x}^{i}(\overline{x}^{\ell+(t-1)i}[t]_{\overline{x}}+1)=0,$$
which is equivalent to $\overline{x}^{t-1}=1$ and $\overline{x}^{\ell}[t]_{\overline{x}}+1=0$. Noting that (\ref{eq:x0}) implies $\overline{x}^{m}=1$, we see that (\ref{eq:cons}) turns to be equivalent to (\ref{eq:x1}) and (\ref{eq:x2}).

If $\sum\limits_{i}^{m}a_{i}\omega_{i}=0$, then $\sum\limits_{i=1}^{m}a_{i}$ is even, hence $\sum\limits_{i}^{m}a_{i}\omega_{i}$ is equal to a linear combination of $2\omega_{m}$ and the $\theta_{i}$'s.
Thus $\Gamma\cong\Gamma_{Q}$ and $\mathcal{M}\cong\mathcal{M}_{Q}$.
\end{proof}

The problem is reduced to determining all ideals $Q\in\mathcal{I}(N,\ell,m,t)$, or, equivalently, determining all quotients $F_{N}(Q)$ with $Q\in\mathcal{I}(N,\ell,m,t)$.
In the following subsections, we separately discusss three cases, and in each case, we identify $F_{N}(Q)$ and $\mathcal{M}_{Q}$,
with more explicit descriptions.

\subsection{$t=m-1$}

First assume $\ell=(t-1,m)=(2,m)$. By (\ref{eq:x0}), (\ref{eq:x1}), (\ref{eq:x2}), $\overline{x}=1$ and then $m=0$ in $F_{N}(Q)$, i.e., $N\mid m$; on the other hand, $[N]_{\overline{x}}=N=0$, hence by (\ref{eq:distinct}), $N=m$.
Thus $F_{N}(Q)\cong\mathbb{Z}_{m}$, and
$$\Gamma_{Q}\cong\langle\mathbb{Z}_{m},\omega\mid 2\omega=-\ell\rangle\cong\mathbb{Z}_{m}\times\mathbb{Z}_{2}, \ \ \ \ \
\overline{x}^{i}\mapsto (1,0), \ \ \omega\mapsto (-\ell/2,1),$$
hence
\begin{align}
\mathcal{M}_{Q}\cong\mathcal{CM}(\mathbb{Z}_{m}\times\mathbb{Z}_{2},\{(i-\ell/2,1)\colon 1\le i\le m\}).
\end{align}

Now assume $2\mid m$ and $\ell=(t-1,m)/2=1$. Then by (\ref{eq:x1}), $\overline{x}^{2}=1$.
Writing $m=2m_{1}$, (\ref{eq:x0}) is equivalent to
$m_{1}(\overline{x}+1)=0.$
Since $[2N]_{\overline{x}}=N(\overline{x}+1)=0$, by (\ref{eq:distinct}) we have $m\mid 2N$, i.e., $m_{1}\mid N$.

If $\overline{x}=a\in\mathbb{Z}_{N}$, then $a^{2}=1$, and there are isomorphisms
\begin{align}
\Gamma_{Q}&\cong\langle\mathbb{Z}_{N},\omega\mid 2\omega=-1\rangle\cong\mathbb{Z}_{2N}, \ \ \ \ \
\overline{x}^{i}\mapsto 2a^{i},  \ \ \omega\mapsto -1,   \nonumber \\
\mathcal{M}_{Q}&\cong\mathcal{CM}(\mathbb{Z}_{2N},\{2[i]_{a}-1\colon 1\le i\le m\}).
\end{align}

If $\overline{x}\notin\mathbb{Z}_{N}$, then as an abelian group,
\begin{align*}
F_{N}(Q)=\langle 1,\overline{x}+1\rangle\cong\mathbb{Z}_{N}\times\mathbb{Z}_{m_{1}}, \ \  \ \
\overline{x}^{2j}\mapsto (1,0),  \ \  \overline{x}^{2j+1}\mapsto (-1,1),
\end{align*}
through which $[i]_{\overline{x}}$ is sent to $(((-1)^{i-1}+1)/2,[i/2])$,
($[c]$ denotes the largest integer not exceeding $c$);
in particular, $[\ell]_{\overline{x}}$ is sent to $(1,0)$.
We have
$$\Gamma_{Q}\cong\langle\mathbb{Z}_{N}\times\mathbb{Z}_{m_{1}},\omega\mid 2\omega=(-1,0)\rangle\cong\mathbb{Z}_{2N}\times\mathbb{Z}_{m_{1}};$$
the later isomorphism sends $(u,v)$ to $(2u,v)$ and sends $\omega$ to $(-1,0)$.
Thus
\begin{align}
\mathcal{M}_{Q}\cong\mathcal{CM}(\mathbb{Z}_{2N}\times\mathbb{Z}_{m_{1}},\{((-1)^{i-1},[i/2])\colon 1\le i\le m\}).
\end{align}

As can be checked, the above results recover Theorem 7.1 of \cite{CJT07-t}, which completely classifies RBCM$_{-1}$'s on abelian groups.

\subsection{$1<t<m-1$ and $\ell=(t-1,m)/2$}

Suppose $N=\prod\limits_{p\in\Delta}p^{k_{p}}$, then
$$F_{N}(Q)\cong\prod\limits_{p\in\Delta}F_{p^{k_{p}}}(Q_{p}), \qquad  \text{with} \quad  Q_{p}=\mathcal{R}_{p^{k_{p}}}(Q).$$
By (\ref{eq:x1}) and Proposition \ref{prop:iso-of-ring} (a), there is a ring isomorphism
$$F_{N}(Q)\cong\prod\limits_{p\in\Delta}\prod\limits_{\lambda\in\Lambda(p,2\ell)}F_{p^{k_{p}}}(Q_{p}(\lambda)).$$

For $p\in\Delta$ and $\lambda\in\Lambda(p,2\ell)$,
consider (\ref{eq:x0}), (\ref{eq:x1}) and (\ref{eq:x2}) in $F_{p^{k_{p}}}(Q_{p}(\lambda))$.

If $\lambda\ne\o$, i.e., $\psi_{\lambda}\neq x-1$, then by (\ref{eq:csq}), $\overline{x}^{\ell}=-1$; conversely, the single condition $\overline{x}^{\ell}=-1$ ensures  (\ref{eq:x0}), (\ref{eq:x1}), (\ref{eq:x2}) to hold.

If $p\neq 2$ and  $\lambda=\o$, then (\ref{eq:x2}) implies $p\mid t+1$, so that $p\nmid t-1$; by Remark \ref{rmk:divide}, $$\mathcal{R}_{p}(Q_{p}(\o))=(x-1),$$
hence $\overline{x}=1+p^{w}\chi\in\mathbb{Z}_{p^{k_{p}}}$ for some $w,\chi$ with $0<w\le k_{p}$ and $p\nmid \chi$.
Since by (\ref{eq:csq}), $\overline{x}-1=0$, we have $p^{w}\in Q_{p}(\o)$.
Then (\ref{eq:x0}) is equivalent to $p^{w}\mid m$ and (\ref{eq:x2}) is equivalent to $p^{w}\mid t+1$.

Consequently, we may rewrite
\begin{align}
F_{N}(Q)\cong F_{N_{1}}(P)\times\mathbb{Z}_{N_{2}}\times F_{2^{k}}(P'),
\end{align}
where:
\begin{itemize}
  \item $k=\exp(F_{2^{k_2}}(Q_{2}(\o)))\le k_2$;
  \item $N_{2}$ is an odd divisor of $(m,t+1)$;
  \item the least common multiple of $2^{k}$, $N_{1}$ and $N_{2}$ is equal to $N$;
  \item $P$ contains $x^{\ell}+1$ and coprime to $x-1$;
  \item $P'=\mathcal{R}_{2^{k}}(Q_{2}(\o))$.
\end{itemize}
Suppose $m=2^{s}m',2\ell=2^{w}\ell'$ with $m',\ell'$ odd, $\ell'\mid m'$ and $w\le s$. It follows from (\ref{eq:x0}) and (\ref{eq:x1}) that
$$0=[m]_{\overline{x}}=\left[\frac{m'}{\ell'}\cdot 2^{s}\ell'\right]_{\overline{x}}=\frac{m'}{\ell'}[2^{s}\ell']_{\overline{x}},$$
so $[2^{s}\ell']_{\overline{x}}=0$ in $F_{2^{k}}(P')$; furthermore, since  $[2\ell]_{\overline{x}}=0$ in $F_{N_{1}}(P)$, we have $[2^{s}\ell'N_{2}]_{\overline{x}}=0$ in $F_{N}(Q)$. Hence by (\ref{eq:distinct}), $m'\mid \ell'N_{2}$, which together with $(\ell',N_{2})=1$ implies
\begin{align}
m'=\ell'N_{2}.    \label{eq:m-ell-N}
\end{align}

If $d_{2}(t-1)\ge s$ so that $w=s$, then in $F_{2^{k}}(P')$, $[2\ell]_{\overline{x}}=0$, hence by (\ref{eq:x2}), $\overline{x}^{\ell}=-1$.
By Example \ref{exmp} and Theorem \ref{thm:ideal},
\begin{align}
P'=M_{v,k-1}^{\o}=(x^{2^{s-1}}+1,2^{k-1}(x-1)^{v}),
\end{align}
for some $v\le 2^{s-1}$.
(Note that when $s=1$, $P'=(x+1)$ so that $F_{2^{k}}(P')\cong\mathbb{Z}_{2^{k}}$).
We can define an isomorphism (of abelian groups)
$$F_{N}(Q)\cong F_{N_{1}}(P)\times\mathbb{Z}_{m/2\ell}\times F_{2^{k}}(M_{v,k-1}^{\o}), \qquad
\overline{x}^{i}\mapsto(\overline{x}^{i}(\overline{x}-1),1,\overline{x}^{i});$$
it sends $[\ell]_{\overline{x}}$ to $(-2,\ell,[\ell]_{\overline{x}})$.
Take
\begin{align}
\beta\in F_{2^{k}}(M_{v,k-1}^{\o}) \qquad \text{with}  \quad 2\beta=[\ell]_{\overline{x}},
\end{align}
then we have an isomorphism
$$\Gamma_{Q}\cong F_{N_{1}}(P)\times\mathbb{Z}_{m/\ell}\times F_{2^{k}}(M_{v,k-1}^{\o}),$$
through which $[i]_{\overline{x}}+\omega$ is sent to $(\overline{x}^{i},2i-\ell,[i]_{\overline{x}}-\beta)$, thus
\begin{align}
\mathcal{M}_{Q}\cong\mathcal{CM}(F_{N_{1}}(P)\times\mathbb{Z}_{m/\ell}\times F_{2^{k}}(M_{v,k-1}^{\o}),
\{(\overline{x}^{i},2i-\ell,[i]_{\overline{x}}-\beta)\colon 1\le i\le m\}).
\end{align}

\medskip

Now assume $d_{2}(t-1)<s$.

Writing $2\ell=b(t-1)+cm$ with $2\nmid b$, we have
$$[\ell]_{\overline{x}}(\overline{x}^{\ell}+1)=[2\ell]_{\overline{x}}
=[b(t-1)]_{\overline{x}}+\overline{x}^{b(t-1)}[cm]_{\overline{x}}=b[t-1]_{\overline{x}}=-b(\overline{x}^{\ell}+1),$$
hence by (\ref{eq:csq}),
$$(\ell+b)(\overline{x}^{\ell}+1)=0.$$
If $2\mid\ell$, then $\overline{x}^{\ell}+1=0$ so that $[2\ell]_{\overline{x}}=0$, but this contradicts (\ref{eq:distinct}). Thus $2\nmid\ell$, and by Remark \ref{rmk:divide}, $\mathcal{R}_{2}(P')=x-1$, i.e., $\overline{x}=a$ for some $a\in\mathbb{Z}_{2^{k}}$.
\begin{enumerate}
  \item[\rm (i)] If $a=1$, then $2^{k}\mid t+1, m$, and (\ref{eq:distinct}) implies $s=k$.
  \item[\rm (ii)] If $k\ge 2$ and $a=-1$, then $[2]_{\overline{x}}=0$, contradicting  (\ref{eq:distinct}).
  \item[\rm (iii)] If $k>2$ and $a=1+2^{w}\chi\neq\pm 1$ with $w\ge 2$ and $2\nmid\chi$, then
        $4\nmid a^{\ell}+1$, hence by (\ref{eq:csq}), $a-1\equiv 0\pmod{2^{k-1}}$. Then (\ref{eq:x0}) implies $s\ge k$, which together with (\ref{eq:distinct}) implies $s=k$;  (\ref{eq:x2}) is equivalent to $2^{k}\mid t+1$.
  \item[\rm (iv)] If $k>2$ and $a=1+2\chi\neq\pm 1$ with $2\nmid\chi$, then $4\nmid a-1$,
        hence by (\ref{eq:csq}) we have $a\equiv -1\pmod{2^{k-1}}$. By (\ref{eq:distinct}), $s=2$.
\end{enumerate}
Let $\tilde{a}\in\{1,\ldots,2^{k}N_{2}\}$ be the unique integer with
\begin{align}
\tilde{a}\equiv 1\pmod{N_{2}} \qquad \text{and} \qquad \tilde{a}\equiv a\pmod{2^{k}}.
\end{align}
We can define an isomorphism (of abelian groups)
\begin{align}
F_{N}(Q)&\cong F_{N_{1}}(P)\times\mathbb{Z}_{2^{k}N_{2}}, \qquad
\overline{x}^{i}\mapsto(\overline{x}^{i}(\overline{x}-1),\tilde{a}^{i}),
\end{align}
by which $[\ell]_{\overline{x}}$ is sent to $(\overline{x}^{\ell}-1,2[\ell]_{\tilde{a}})=(-2,[\ell]_{\tilde{a}})$.
Then
$$\Gamma_{Q}\cong\langle F_{N_{1}}(P)\times\mathbb{Z}_{2^{k}N_{2}},\omega\mid 2\omega=(2,-[\ell]_{\tilde{a}})\rangle\cong F_{N_{1}}(P)\times\mathbb{Z}_{2^{k+1}N_{2}},$$
where the later isomorphism sends $(\overline{\alpha},u)$ to $(\overline{\alpha},2u)$ and sends $\omega$ to $(1,[\ell]_{\tilde{a}})$.
Thus, (recalling (\ref{eq:m-ell-N}),) $\mathcal{M}_{Q}$ is isomorphic to
\begin{align}
\mathcal{M}(P;k,a):=\mathcal{CM}(F_{N_{1}}(P)\times\mathbb{Z}_{2^{k+1}N_{2}},\{(\overline{x}^{i}, 2[i]_{\tilde{a}}-[\ell]_{\tilde{a}})\colon 1\le i\le m\}).
\end{align}

\begin{thm}
Let $\mathcal{M}$ be an $m$-valent type I RBCM$_{t}$ on an abelian group, then exactly one of the following cases occurs:
\begin{enumerate}
  \item[{\rm(a)}] $d_{2}(t-1)\ge d_{2}(m)$, and
                  $$\mathcal{M}\cong\mathcal{CM}(F_{N_{1}}(P)\times\mathbb{Z}_{m/\ell}\times F_{2^{k}}(M_{v,k-1}^{\o}),
                  \{(\overline{x}^{i},2i-\ell,[i]_{\overline{x}}-\beta)\colon 1\le i\le m\})$$
                  for a unique triple $(P,k,v)$ with $P\in\mathcal{I}(N_{1},\ell)$;
  \item[{\rm(b)}] $d_{2}(t-1)=1, d_{2}(m)=2$, and $\mathcal{M}\cong\mathcal{M}(P;k,2^{k-1}-1)$ with $k>2$;
  \item[{\rm(c)}] $d_{2}(t+1)\ge d_{2}(m)=k\ge 2$, and $\mathcal{M}\cong\mathcal{M}(P;k,1)$;
  \item[{\rm(d)}] $d_{2}(t+1)\ge d_{2}(m)=k>2$, and $\mathcal{M}\cong\mathcal{M}(P;k,2^{k-1}+1)$;
\end{enumerate}
in {\rm(b)}, {\rm(c)}, {\rm(d)}, $P\in\mathcal{I}(N_{1},\ell)$ is a unique ideal coprime to $x-1$.
\end{thm}

For abelian 2-groups, by (\ref{eq:m-ell-N}), $\ell'=m'$; the assumption $t\not\equiv\pm 1\pmod{m}$ implies $d_{2}(t-1)<d_{2}(m)$;
suppose $N_1=2^{k'}$.
\begin{cor}
Let $\mathcal{M}$ be an $m$-valent type I RBCM$_{t}$ on an abelian 2-group,
then one of the following two cases occurs:
\begin{itemize}
  \item $d_{2}(t+1)\ge d_{2}(m)=k\ge 2$, and
        $$\mathcal{M}\cong\mathcal{CM}(F_{2^{k'}}(P)\times\mathbb{Z}_{2^{k+1}},\{(\overline{x}^{i}, 2i-\ell)\colon 1\le i\le m\});$$
  \item $d_{2}(t+1)\ge d_{2}(m)=k>2$, and
        $$\mathcal{M}\cong\mathcal{CM}(F_{2^{k'}}(P)\times\mathbb{Z}_{2^{k+1}},\{(\overline{x}^{i}, 2[i]_{a}-[\ell]_{a})\colon 1\le i\le m\});$$
\end{itemize}
in both cases $P\in\mathcal{I}(2^{k'},\ell)$ is a unique ideal coprime to $x-1$, and in the second case, $a=2^{k-1}+1$.
\end{cor}

\subsection{$1<t<m-1$ and $\ell=(t-1,m)$}

Write $\ell=a(t-1)+bm$, then $[\ell]_{\overline{x}}=a[t-1]_{\overline{x}}+b[m]_{\overline{x}}=-2a$, and
$$\left[\frac{N}{2}\ell\right]_{\overline{x}}=\frac{N}{2}[\ell]_{\overline{x}}=-Na=0,$$
hence by (\ref{eq:distinct}), $m\mid N\ell/2$.
On the other hand,
$$-2\frac{m}{\ell}=\frac{m}{\ell}[t-1]_{\overline{x}}=\left[\frac{m}{\ell}(t-1)\right]_{\overline{x}}=\frac{t-1}{\ell}[m]_{\overline{x}}=0,$$
hence $N\mid 2m/\ell$.
Thus
\begin{align}
2m=N\ell.    \label{eq:m-N-ell}
\end{align}
This implies $2\mid N$.
By (\ref{eq:x1}) and (\ref{eq:csq}),
\begin{align}
2(\overline{x}-1)=0, \hspace{5mm}  {\rm i.e.}, \ \ \  2\overline{x}=2.  \label{eq:condition3}
\end{align}

Suppose $N=2^{k}N'$ with $k\ge 1$ and $2\nmid N'$, then
\begin{align}
F_{N}(Q)\cong F_{N'}(Q_{N'})\times F_{2^{k}}(Q_{2}), \qquad \text{with}  \quad   Q_{n}=\mathcal{R}_{n}(Q).
\end{align}

By (\ref{eq:condition3}), $\overline{x}=1$ in $F_{N'}(Q_{N'})$ so that $F_{N'}(Q_{N'})\cong\mathbb{Z}_{N'}$.

If $d_{2}(t-1)\ge d_{2}(m)$, then in $F_{2^{k}}(Q_{2})$, $[t-1]_{\overline{x}}=0$, hence by (\ref{eq:x2}), $2=0$, i.e., $k=1$;
conversely, if $k=1$, then (\ref{eq:m-N-ell}) implies $d_{2}(t-1)\ge d_{2}(m)$. Under the assumption $d_{2}(t-1)\ge d_{2}(m)$ or $k=1$, (\ref{eq:x0}), (\ref{eq:x1}), (\ref{eq:x2}) are equivalent to the single condition $[\ell]_{x}\in Q_{2}$, which, combining with (\ref{eq:distinct}), is in turn equivalent to $(x-1)Q_{2}\in\mathcal{I}(2,\ell)$.
We have isomorphisms
$$\Gamma_{Q}\cong\langle\mathbb{Z}_{N'}\times F_{2}(Q_{2}),\omega\mid 2\omega=-(\ell,0)\rangle
\cong\mathbb{Z}_{N}\times F_{2}(Q_{2});$$
the later sends $(u,\overline{\alpha})$ to $(2u,\overline{\alpha})$ and sends $\omega$ to $(-\ell,0)$,
thus
\begin{align}
\mathcal{M}_{Q}\cong\mathcal{CM}(\mathbb{Z}_{N}\times F_{2}(Q_{2}),\{(2i-\ell,[i]_{\overline{x}})\colon 1\le i\le m\}).
\end{align}

\bigskip

Now assume $s=d_{2}(t-1)<d_{2}(m)$ and $k>1$.
By Proposition \ref{prop:iso-of-ring},
$$F_{2^{k}}(Q_{2})\cong\prod\limits_{\lambda\in\Lambda(2,\ell)}F_{2^{k}}(Q_{2}(\lambda)).$$
If $\lambda\neq \o$, then by (\ref{eq:condition3}), $2\in Q_{2}(\lambda)$, hence $F_{2^{k}}(Q_{2}(\lambda))\cong F_{2}((\vartheta_{\lambda}))$
for some $\vartheta_{\lambda}\in\mathbb{Z}_{2}[x]$ coprime to $x-1$, and (\ref{eq:x0}), (\ref{eq:x1}), (\ref{eq:x2}) are equivalent to the single condition $\vartheta_{\lambda}\mid [\ell]_{x}$. Let $\vartheta=\prod\limits_{\lambda\neq\o}\vartheta_{\lambda}$.
Then
\begin{align}
F_{2^{k}}(Q_{2})\cong F_{2}((\vartheta))\times F_{2^{k}}(Q_{2}(\o)).
\end{align}

Suppose $\mathcal{R}_{2}(Q_{2}(\o))=(x-1)^{u}$ with $u>0$.

If $u>1$, then $k=2$. This because, by (\ref{eq:x2}),
$\mathcal{R}_{2}([t]_{x}+1)=[t-1]_{x}$ can be divided by $(x-1)^{u}$,
hence
$$(x-1)^{u+1}\mid x^{t-1}-1 \qquad   \text{in }  \quad \mathbb{Z}_{2}[x],$$
implying $u+1\le 2^{s}$ by Remark \ref{rmk:divide}, so $s>1$ and $4\nmid t+1$; on the other hand, by (\ref{eq:x2}) and (\ref{eq:condition3}),
$0=2(1+[t]_{\overline{x}})=2(t+1),$ hence
$$k-1\le d_{2}(t+1)=1.$$
Take $\alpha\in Q_{2}(\o)$ with $\mathcal{R}_{2}(\alpha)=(x-1)^{u}$, then each element of $Q_{2}(\o)$ can be written in the form $\alpha\beta+2\gamma$. Let $2a$ (with $a\in\{0,1\}$) be the remainder of $2\alpha$ divided by $2(x-1)$, then,
noting (\ref{eq:condition3}), we see that actually
$$Q_{2}(\o)=((x-1)^{u}+2a,2(x-1)).$$
Also by (\ref{eq:condition3}),
$$[2^{s}]_{\overline{x}}=\sum\limits_{i=0}^{2^{s}-1}\overline{x}^{i}=\sum\limits_{i=0}^{2^{s}-1}\sum\limits_{j=0}^{i}{i\choose j}(\overline{x}-1)^{j}=\sum\limits_{j=1}^{2^{s}}{2^{s}\choose j}(\overline{x}-1)^{j-1}=(\overline{x}-1)^{2^{s}-1}.$$
Since by (\ref{eq:x2}) $[t-1]_{\overline{x}}\ne 0$, we have
$(\overline{x}-1)^{2^{s}-1}=[2^{s}]_{\overline{x}}\ne 0$. Consequently, $a=1$; moreover, $u=2^{s}-1$: otherwise $(\overline{x}-1)^{2^{s}-1}$
would vanish, as
$$(\overline{x}-1)^{2^{s}-1}=(\overline{x}-1)^{2^{s}-1-u}((\overline{x}-1)^{u}+2)-2(\overline{x}-1)^{2^{s}-1-u}.$$
Thus $Q_{2}(\o)=P_{(s)}$, with
\begin{align}
P_{(s)}=((x-1)^{2^{s}-1}+2,2(x-1)).
\end{align}
One can check that (\ref{eq:x0}), (\ref{eq:x1}) and (\ref{eq:x2}) all hold in $F_{4}(P_{(s)})$;
note that
$$[\ell]_{\overline{x}}=\frac{\ell}{2^{s}}\cdot[2^{s}]_{\overline{x}}=\frac{\ell}{2^{s}}\cdot 2=2.$$
Now we have isomorphisms
\begin{align*}
\Gamma_{Q}&\cong\langle\mathbb{Z}_{N'}\times F_{2}((\vartheta))\times F_{4}(P_{(s)}),\omega\mid 2\omega=-(\ell,0,2)\rangle \\
&\cong\mathbb{Z}_{N'}\times\mathbb{Z}_{2}\times F_{2}((\vartheta))\times F_{4}(P_{(s)});
\end{align*}
the later is defined as to send $(u,\alpha,\beta)$ to $(u,0,(\overline{x}-1)\alpha,\overline{x}^{-1}\beta)$ and send $\omega$ to $(-\ell/2,1,0,-\overline{x}^{-1})$, hence the composite sends $[i]_{\overline{x}}+\omega$ to $(i-\ell/2,1,\overline{x}^{i},[i-1]_{\overline{x}})$. Thus $\mathcal{M}_{Q}$ is isomorphic to
\begin{align}
\mathcal{CM}(\mathbb{Z}_{N'}\times\mathbb{Z}_{2}\times F_{2}((\vartheta))\times F_{4}(P_{(s)}),
\{(i-\ell/2,1,\overline{x}^{i},[i-1]_{\overline{x}})\colon 1\le i\le m\}).
\end{align}

\begin{rmk}  \label{rmk:P(s)}
\rm As an abelian group, $F_{4}(P_{(s)})$ is generated by $(\overline{x}-1)^{i},i=0,\ldots,2^{s}-2$, with relations $4=0$ and $2(\overline{x}-1)^{i}=0, i=1,\ldots, 2^{s}-2$, so
\begin{align}
F_{4}(P_{(s)})\cong\mathbb{Z}_{4}\times\mathbb{Z}_{2}^{2^{s}-2}.
\end{align}
\end{rmk}

\medskip

If $u=1$, then in $F_{2^{k}}(Q_{2}(\o))$, $\overline{x}=a\in\mathbb{Z}_{2^{k}}$ with $2(a-1)=0$, i.e., $a\in\{1,1+2^{k-1}\}$.
Now (\ref{eq:x0}) implies $d_{2}(m)\ge k$; on the other hand,
$$[2^{k}]_{\overline{x}}=2^{k-1}(1+a)=0,$$
hence by (\ref{eq:distinct}), $d_{2}(m)\le k$. Thus $d_{2}(m)=k$, and by (\ref{eq:m-N-ell}), $s=1$.
By (\ref{eq:x2}),
$$d_{2}(t+1)\left\{\begin{array}{ll} \ge k, & a=1, \\ =k-1\ (\text{\ so\ that\ } k>2), &a=1+2^{k-1} \end{array}\right.$$
When $d_{2}(t+1)\ge k$, we have isomorphisms
\begin{align*}
F_{N}(Q)\cong \mathbb{Z}_{N'}\times F_{2}((\vartheta))\times\mathbb{Z}_{2^{k}}\cong\mathbb{Z}_{N}\times F_{2}((\vartheta)),
\end{align*}
through which $\overline{x}^{i}$ is sent to $(1,\overline{x}^{i}(\overline{x}-1))$, and
\begin{align*}
\Gamma_{Q}\cong\langle\mathbb{Z}_{N}\times F_{2}((\vartheta)),\omega\mid 2\omega=-(\ell,0)\rangle
\cong\mathbb{Z}_{N}\times\mathbb{Z}_{2}\times F_{2}((\vartheta)),
\end{align*}
the later sending $(u,\alpha)$ to $(u,0,\alpha)$ and sending $\omega$ to $(-\ell/2,1,0)$, thus
\begin{align}
\mathcal{M}_{Q}\cong
\mathcal{CM}(\mathbb{Z}_{N}\times\mathbb{Z}_{2}\times F_{2}((\vartheta)),\{(i-\ell/2,1,\overline{x}^{i})\colon 1\le i\le m\}).
\end{align}
When $d_{2}(t+1)=k-1$, let
\begin{align}
\tilde{a}=\frac{N}{2}+1,
\end{align}
then similarly as above,
\begin{align}
\mathcal{M}_{Q}&\cong\mathcal{CM}(\mathbb{Z}_{N}\times\mathbb{Z}_{2}\times F_{2}((\vartheta)),
\{([i]_{\tilde{a}}-[\ell]_{\tilde{a}}/2,1,\overline{x}^{i})\colon 1\le i\le m\}).
\end{align}

\begin{thm} \label{thm:II-RBCMt}
Let $\mathcal{M}$ be an $m$-valent type II RBCM$_{t}$ on an abelian group of exponent $N$, then $N=2m/\ell$ and one of the following cases occurs:
\begin{enumerate}
  \item[\rm(a)] $d_{2}(t-1)\ge d_{2}(m)$, and $$\mathcal{M}\cong\mathcal{CM}(\mathbb{Z}_{2m/\ell}\times F_{2}((\vartheta)),\{(2i-\ell,[i]_{\overline{x}})\colon 1\le i\le m\})$$
      for a unique $\vartheta\in\mathbb{Z}_{2}[x]$ such that $((x-1)\vartheta)\in\mathcal{I}(2,\ell)$.
  \item[\rm(b)] $d_{2}(m)-1=d_{2}(t-1)=s>1$, and
      $$\mathcal{M}\cong\mathcal{CM}(\mathbb{Z}_{m/2\ell}\times\mathbb{Z}_{2}\times F_{2}((\vartheta))\times F_{4}(P_{(s)}),
      \{(i-\ell/2,1,\overline{x}^{i},[i-1]_{\overline{x}})\colon 1\le i\le m\});$$
  \item[\rm(c)] $d_{2}(t+1)\ge d_{2}(m)\ge 2$, and
      $$\mathcal{M}\cong
      \mathcal{CM}(\mathbb{Z}_{2m/\ell}\times\mathbb{Z}_{2}\times F_{2}((\vartheta)),\{(i-\ell/2,1,\overline{x}^{i})\colon 1\le i\le m\});$$
  \item[\rm(d)] $d_{2}(t+1)+1=d_{2}(m)>2$, and
      $$\mathcal{M}\cong
      \mathcal{CM}(\mathbb{Z}_{2m/\ell}\times\mathbb{Z}_{2}\times F_{2}((\vartheta)),
      \{([i]_{\tilde{a}}-[\ell]_{\tilde{a}}/2,1,\overline{x}^{i})\colon 1\le i\le m\});$$
\end{enumerate}
in {\rm(b)}, {\rm(c)}, {\rm(d)}, $\vartheta$ is a unique polynomial with $x-1\nmid \vartheta$ and $(\vartheta)\in\mathcal{I}(2,\ell)$.
\end{thm}

Since the ideals in $\mathcal{I}(2,\ell)$ are all known, we can say that type II RBCM$_{t}$'s on abelian groups are completely classified.

\begin{exmp}
\rm We can recover the results of Theorem 6.7 of \cite{CJT07} by considering some special cases. Let $m=\tilde{t}^2-1$, let $\mathcal{M}^{\pm}$ denote the RBCM$_{\pm\tilde{t}}$ on $\mathbb{Z}_{2(\tilde{t}\pm 1)}\times \mathbb{Z}_2^{\tilde{t}-1\mp 1}$ constructed there, and let $\mathcal{M}^a(\vartheta)$ (resp. $\mathcal{M}^b(\vartheta)$, $\mathcal{M}^d(\vartheta)$) denote the RBCM$_t$ given in the case (a) (resp. (b), (d)) of Theorem \ref{thm:II-RBCMt}.

First suppose $\tilde{t}$ is even so that $m$ is odd.
\begin{itemize}
  \item Let $t=\tilde{t}$, then $\ell=\tilde{t}-1, m/\ell=\tilde{t}+1$ and $d_2(t-1)\ge d_2(m)$. Let $b_1,\ldots,b_{t-2}$ denote the generators of $\mathbb{Z}_{2(t+1)}\times\mathbb{Z}_2^{t-3}$ as in \cite{CJT07} Lemma 6.6. There is an isomorphism
      $$\mathbb{Z}_{2(t+1)}\times\mathbb{Z}_2^{t-3}\to\mathbb{Z}_{t+1}\times F_2(([t-1]_x))$$
      sending $b_1$ to $(1,1)$ and sending $b_j$ to $(0,\overline{x}^{j-1}+\overline{x}^{j-2})$ for $j=2,\ldots,t-2$.
      Then it can be verified that
      $\mathcal{M}^+\cong\mathcal{M}^a([\ell]_x)$.
  \item Let $t=m-\tilde{t}$, then $\ell=\tilde{t}+1, m/\ell=\tilde{t}-1$ and $d_2(t-1)\ge d_2(m)$. Similarly as above,
      $\mathcal{M}^-\cong\mathcal{M}^a([\ell]_x)$.
\end{itemize}

Second, suppose $\tilde{t}$ is odd with $4\nmid \tilde{t}-1$.
\begin{itemize}
  \item Let $t=\tilde{t}$, then $\ell=\tilde{t}-1$, $d_{2}(t+1)+1=d_{2}(m)>2$, and
      $\mathcal{M}^+\cong\mathcal{M}^d(([\ell/2]_x)^2)$.
  \item Let $t=m-\tilde{t}$, then $\ell=\tilde{t}+1$, $d_{2}(m)-1=d_{2}(t-1)=s>1$, and
      $\mathcal{M}^-\cong\mathcal{M}^b(([\ell/2^s]_x)^{2^s})$;
      in this case we need to verify that
      $$\mathbb{Z}_{2(\tilde{t}-1)}\times\mathbb{Z}_2^{\tilde{t}-1}
      \cong\mathbb{Z}_{(\tilde{t}-1)/2}\times F_2((\vartheta))\times F_4(P_{(s)})
      \cong\mathbb{Z}_{(\tilde{t}-1)/2}\times F_2((\vartheta))\times\mathbb{Z}_4\times\mathbb{Z}_2^{2^s-2},$$
      with $\vartheta=([\ell/2^s]_x)^{2^s}$, which is a routine job.
\end{itemize}

Finally, suppose $\tilde{t}$ is odd with $4\mid \tilde{t}-1$.
\begin{itemize}
  \item Let $t=\tilde{t}$, then $\ell=\tilde{t}-1$, $d_{2}(m)-1=d_{2}(t-1)=s>1$, and
      $\mathcal{M}^+\cong\mathcal{M}^b(([\ell/2^s]_x)^{2^s})$.
  \item Let $t=m-\tilde{t}$, then $\ell=\tilde{t}-1$, $d_{2}(t+1)+1=d_{2}(m)>2$, and
      $\mathcal{M}^-\cong\mathcal{M}^d(([\ell/2]_x)^{2})$.
\end{itemize}
\end{exmp}

For abelian 2-groups, by (\ref{eq:m-N-ell}) and the assumption $t\not\equiv\pm 1\pmod{m}$  we have $d_{2}(t-1)<d_{2}(m)$ and $d_{2}(N)\ge 2$.

\begin{cor}
Let $\mathcal{M}$ be an $m$-valent type II RBCM$_{t}$ on an abelian 2-group of exponent $2^{k}$,
then $m=2^{k-1}\ell$ and one of the following cases occurs:
\begin{itemize}
  \item $k=2$, $d_{2}(m)-1=d_{2}(t-1)=s>1$, and
       $$\mathcal{M}\cong\mathcal{CM}(F_{4}(P_{(s)})\times\mathbb{Z}_{2}\times F_{2}((\vartheta)),
       \{([i-1]_{\overline{x}},1,\overline{x}^{i})\colon 1\le i\le m\});$$
  \item $d_{2}(t+1)\ge d_{2}(m)=k\ge 2$, and
      $$\mathcal{M}\cong
      \mathcal{CM}(\mathbb{Z}_{2^{k}}\times\mathbb{Z}_{2}\times F_{2}((\vartheta)),\{(i-\ell/2,1,\overline{x}^{i})\colon 1\le i\le m\});$$
  \item $d_{2}(t+1)+1=d_{2}(m)=k>2$, and
      $$\mathcal{M}\cong
      \mathcal{CM}(\mathbb{Z}_{2^{k}}\times\mathbb{Z}_{2}\times F_{2}((\vartheta)),
      \{([i]_{\tilde{a}}-[\ell]_{\tilde{a}}/2,1,\overline{x}^{i})\colon 1\le i\le m\});$$
\end{itemize}
here $\vartheta$ is a unique polynomial with $x-1\nmid\vartheta$ and $(\vartheta)\in\mathcal{I}(2,\ell)$.
\end{cor}

\end{document}